\documentclass{amsart}
\usepackage{a4wide,amssymb,color,hyperref}
\usepackage{graphicx} 

\theoremstyle{plain}

\theoremstyle{definition}
\newtheorem{definition}{Definition}[section]

\theoremstyle{remark}

\newtheorem{example}{Example}[section]

\title{Enumerations of 1-rotational Steiner systems}
\author{Ivan Hetman, Taras Banakh, Alex Ravsky}
\date{May 2025}

\begin{document}

\begin{abstract} In this paper new $1$-rotational Steiner systems for different admissible $v,k$ pairs are introduced. In particular, 1-rotational unitals of order $4$ are enumerated.
\end{abstract}

\maketitle

\section{Introduction}

In \cite{Het} results of an algorithm for finding cyclic difference families (CDF) were presented. Since then this algorithm was generalized from cyclic to commutative groups, then to any groups and in the end to sets endowed with group actions. Then while rereading VI.16 chapter of \cite{HoCD}, it occured that cyclic/commutative/any group algorithms can be very easily generalized to the case of so-called 1-rotational designs. After that corresponding changes were applied and this paper will introduce new Steiner designs generated by various groups.

$S(2,5,65)$ 1-rotational unitals of order $4$ for groups of order $64$ intersect with results in \cite{Krc} and \cite{NagyWeb}, still there are lots of missing unitals in \cite{Krc}. Still, due to lots of designs and relatively slow algorithms checking the isomorphism of designs, we didn't check for intersection with \cite{NagyWeb}, but we assume that there might be designs missing there too.

Due to the fact that we have different specifically optimized algorithms for cyclic/commutative/any group cases, some cases will be presented in form multiplier-fingerprint-difference family with common group element form, but results generated by generalized group algorithm are in form fingerprint-difference family. The Cayley tables of this search were generated by the procedure {\tt CayleyTable(IntoLoop(SmallGroup(v - 1,idx)))} in the LOOPS package \cite{LOOPS} of GAP \cite{GAP4}. As most groups are non-commutative, group action ob blocks is applied from the left.

For $1$-rotational difference families the following equivalence relation \cite{HoCD} was applied.
\begin{definition} Two difference families $D=\{B_1,B_2,...,B_m\}$ and $D'=\{B'_1,B'_2,...,B'_m\}$ over a group $G$ are {\em equivalent} if there is an isomorphism $\alpha$ of $G$ such that $B'_i$ is a translate of $\alpha(B_i)$ for all $i$.
\end{definition}
Difference family equivalence is coarser than design isomorphism because equivalent difference families lead to isomorphic Steiner systems, but the converse is not true in general. Examples 2.4, 3.3 and table for unitals of order $4$ show that sometimes it filters well, but sometimes isomorphic Steiner systems are not filtered. We are filtering designs using difference families equivalence relation.

\section{Enumerations of 1-rotational designs for k = 3}

While cyclic difference families (CDF) start occuring immediately, 1-rotational analogue for cyclic groups start occuring later and with "holes". Some results will be in the form "there is no 1-rotational designs" for specific $v,k$ pairs. At some point we will stop enumerations, still for small $k \le 5$ the algorithm still has capacity to generate full enumeration for next admissible $v$. Sometimes a group that is acting appears in prefix.

\begin{center}
\begin{tabular}{||c c c c||} 
 \hline
 GAP ID & Structure & \# of diff families & Comments \\ [0.5ex] 
 \hline\hline
 {\tt SmallGroup(8,1)} & {\tt C8} & 1 & - \\ 
 {\tt SmallGroup(8,4)} & {\tt Q8} & 1 & - \\ 
 \hline
 {\tt SmallGroup(18,5)} & {\tt C6 x C3} & 1 & - \\
 \hline
 {\tt SmallGroup(24,3)} & {\tt SL(2,3)} & 7 & - \\
 \hline
 {\tt SmallGroup(26,2)} & $\mathbb Z_{26}$ & 36 & - \\
 \hline
 {\tt SmallGroup(32,1)} & {\tt C32} & 242 & - \\
 {\tt SmallGroup(32,20)} & {\tt Q32} & 84 & - \\
 \hline
 {\tt SmallGroup(42,2)} & {\tt C2 x (C7 : C3)} & 14688 & - \\
 \hline
 {\tt SmallGroup(50,2)} & $\mathbb Z_{50}$ & 550528 & - \\
 {\tt SmallGroup(50,4)} & $\mathbb Z_2 \times \mathbb Z_5 \times \mathbb Z_5$ & 24888 & $84$ with multiplier $3$ \\
 \hline
 {\tt SmallGroup(57,1)} & {\tt C7 : C8} & $\ge 1$ & - \\
 {\tt SmallGroup(57,2)} & {\tt C56} & $\ge 1$ & - \\
 {\tt SmallGroup(57,3)} & {\tt C7 : Q8} & $\ge 1$ & - \\
 {\tt SmallGroup(57,10)} & {\tt C7 x Q8} & $\ge 1$ & - \\
 \hline
 {\tt SmallGroup(72,12)} & {\tt C3 x (C3 : C8)} & $\ge 1$ & - \\
 {\tt SmallGroup(72,13)} & {\tt (C3 x C3) : C8} & $\ge 1$ & - \\
 {\tt SmallGroup(72,14)} & {\tt C24 x C3} & $\ge 1$ & - \\
 {\tt SmallGroup(72,19)} & {\tt (C3 x C3) : C8} & $\ge 1$ & - \\
 {\tt SmallGroup(72,24)} & {\tt (C3 x C3) : Q8} & $\ge 1$ & - \\
 {\tt SmallGroup(72,25)} & {\tt C3 x SL(2,3)} & $\ge 1$ & - \\
 {\tt SmallGroup(72,26)} & {\tt C3 x (C3 : Q8)} & $\ge 1$ & - \\
 {\tt SmallGroup(72,31)} & {\tt (C3 x C3) : Q8} & $\ge 1$ & - \\
 {\tt SmallGroup(72,38)} & {\tt C3 x C3 x Q8} & $\ge 1$ & - \\
 \hline
\end{tabular}
\end{center}

\begin{example} In \cite{Bur} it was conjectured that 1-rotational designs could exist only for $v\equiv 1,3,9,19 \mod 24$. It was confirmed by calculations that there are no 1-rotational designs for $(v,k)\in\{(7,3),(13,3),(15,3),(21,3),(31,3),(37,3),(39,3),(45,3),(55,3),(61,3),(63,3),(67,3),(69,3)\}$.
\end{example}

\begin{example} Two groups generate the same affine plane of order $3$ for $(v,k)=(9,3)$.
\begin{enumerate}
    \item {\tt SmallGroup(8,1) = C8} \{1=432\} [[0, 1, 4], [0, 3, $\infty$]]
    \item {\tt SmallGroup(8,4) = Q8} \{1=432\} [[0, 1, 2], [0, 3, $\infty$]]
\end{enumerate}
\end{example}

\begin{example} There is one 1-rotational design for $(v,k)=(19,3)$ for {\tt SmallGroup(18,5) = C6 x C3}.
\begin{enumerate}
    \item \{0=432, 1=5040\} [[0, 1, $\infty$], [0, 2, 5], [0, 3, 8], [0, 4, 12], [0, 7, 16]]
\end{enumerate}
\end{example}

\begin{example} There are seven 1-rotational designs for $(v,k)=(25,3)$, all for the group \newline {\tt SmallGroup(24,3) = SL(2,3)}.
\begin{enumerate}
\item \{0=192, 1=13008\} [[0, 1, 2], [0, 3, 13], [0, 4, $\infty$], [0, 7, 20], [0, 8, 18], [0, 15, 21]]
\item {\bf \{0=1536, 1=11664\} [[0, 1, 5], [0, 2, 3], [0, 4, $\infty$], [0, 6, 12], [0, 7, 20], [0, 8, 19]]}
\item \{0=1536, 1=11664\} [[0, 1, 5], [0, 2, 3], [0, 4, $\infty$], [0, 6, 14], [0, 12, 16], [0, 15, 21]]
\item \{0=960, 1=12240\} [[0, 1, 2], [0, 3, 16], [0, 4, $\infty$], [0, 6, 23], [0, 8, 22], [0, 15, 21]]
\item \{0=1536, 1=11664\} [[0, 1, 5], [0, 2, 3], [0, 4, $\infty$], [0, 6, 20], [0, 8, 12], [0, 15, 21]]
\item \{0=960, 1=12240\} [[0, 1, 5], [0, 2, 6], [0, 3, 17], [0, 4, $\infty$], [0, 7, 20], [0, 8, 22]]
\item \{0=192, 1=13008\} [[0, 1, 2], [0, 3, 19], [0, 4, $\infty$], [0, 6, 23], [0, 7, 20], [0, 8, 15]]
\end{enumerate}
\end{example}
This contradicts a theorem proved in \cite{Bur} where it is claimed that this design is unique. Computer calculations show that all $7$ obtained designs are non-isomorphic (in spite of the fact that some of them have the same fingerprints). Buratti's design is marked bold.

After this we will stop adding more designs as because for an admissible $v$ and $k=3$ the number of designs grows very quickly. Still, anyone who is interested, can adapt our algorithm and generate full enumerations for bigger $v$.

\section{Enumerations of 1-rotational designs for k = 4}

\begin{center}
\begin{tabular}{||c c c c||} 
 \hline
 GAP ID & Structure & \# of diff families & Comments \\ [0.5ex] 
 \hline\hline
 {\tt SmallGroup(15,1)} & $\mathbb Z_{15}$ & 1 & - \\
 \hline
 {\tt SmallGroup(24,10)} & {\tt C3 x D8} & 1 & - \\
 {\tt SmallGroup(24,12)} & {\tt S4} & 2 & - \\
 \hline
 {\tt SmallGroup(27,3)} & {\tt (C3 x C3) : C3} & 2 & - \\
 {\tt SmallGroup(27,4)} & {\tt C9 : C3} & 1 & - \\
 {\tt SmallGroup(27,5)} & {\tt C3 x C3 x C3} & 1 & - \\
 \hline
 {\tt SmallGroup(39,1)} & {\tt C13 : C3} & 28 & - \\
 {\tt SmallGroup(39,2)} & {\tt C39} & 2 & - \\
 \hline
 {\tt SmallGroup(48,3)} & {\tt (C4 x C4) : C3} & 105 & - \\
 {\tt SmallGroup(48,10)} & {\tt (C3 : C8) : C2} & 16 & - \\
 {\tt SmallGroup(48,11)} & {\tt C4 x (C3 : C4)} & 40 & - \\
 {\tt SmallGroup(48,20)} & {\tt C12 x C4} & 4 & - \\
 {\tt SmallGroup(48,24)} & {\tt C3 x (C8 : C2)} & 26 & - \\
 \hline
 {\tt SmallGroup(51,1)} & $\mathbb Z_{51}$ & 166 & - \\
 \hline
 {\tt SmallGroup(63,1)} & {\tt C7 : C9} & 5508 & - \\
 {\tt SmallGroup(63,2)} & {\tt C63} & 18532 & - \\
 {\tt SmallGroup(63,3)} & {\tt C3 x (C7 : C3)} & 77507 & - \\
 {\tt SmallGroup(63,4)} & {\tt C21 x C3} & 9766 & - \\
 \hline
 {\tt SmallGroup(72,10)} & {\tt C9 x D8} & 38600 & - \\
 {\tt SmallGroup(72,15)} & {\tt ((C2 x C2) : C9) : C2} & 4 & - \\
 {\tt SmallGroup(72,22)} & {\tt (C6 x S3) : C2} & 96 & - \\
 {\tt SmallGroup(72,28)} & {\tt C3 x D24} & 544 & - \\
 {\tt SmallGroup(72,30)} & {\tt C3 x ((C6 x C2) : C2)} & 3308 & - \\
 {\tt SmallGroup(72,37)} & {\tt C3 x C3 x D8} & 30840 & - \\
 {\tt SmallGroup(72,40)} & {\tt (S3 x S3) : C2} & 128 & - \\
 {\tt SmallGroup(72,42)} & {\tt C3 x S4} & 52259 & - \\
 {\tt SmallGroup(72,43)} & {\tt (C3 x A4) : C2} & 32 & - \\
 \hline
\end{tabular}
\end{center}

\begin{example} There are no 1-rotational designs for $(v,k)\in\{(13,4),(37,4),(61,4)\}$.
\end{example}

\begin{example} There is only one 1-rotational design for $(v,k)=(16,4)$ which is obtained by group $\mathbb Z_{15}$: [[0, 1, 3, 7], [0, 5, 10, $\infty$]].
\end{example}

\begin{example} There is only one 1-rotational design up to isomorphism for $(v,k)=(25,4)$ that can be obtained in three ways.
\begin{enumerate}
\item {\tt SmallGroup(24,10) = C3 x D8} \{1=12096, 2=13104\} \newline [[0, 1, 4, 7], [0, 2, 6, 20], [0, 3, 10, $\infty$], [0, 5, 11, 16]]
\item {\tt SmallGroup(24,12) = S4} \{1=12096, 2=13104\} \newline  [[0, 1, 11, 15], [0, 2, 6, 17], [0, 3, 12, 19], [0, 4, 5, 14], [0, 9, 22, $\infty$]]
\item {\tt SmallGroup(24,12) = S4} \{1=12096, 2=13104\} \newline  [[0, 1, 9, 13], [0, 2, 7, 19], [0, 3, 4, 11], [0, 5, 10, 20], [0, 17, 18, $\infty$]]
\end{enumerate}
\end{example}

\begin{example} There are $4$ 1-rotational designs (unitals of order $3$) for $(v,k)=(28,4)$.
\begin{enumerate}
\item {\tt SmallGroup(27,3) = (C3 x C3) : C3} \{0=810, 1=12204, 2=23274\} \newline  [[0, 1, 2, 3], [0, 5, 8, 26], [0, 18, 20, $\infty$]]
\item {\tt SmallGroup(27,3) = (C3 x C3) : C3} \{2=36288\} \newline  [[0, 1, 2, 18], [0, 3, 9, $\infty$], [0, 5, 14, 15]]
\item {\tt SmallGroup(27,4) = C9 : C3} \{1=15120, 2=21168\} \newline  [[0, 1, 2, 20], [0, 3, 14, 26], [0, 8, 22, $\infty$]]
\item {\tt SmallGroup(27,5) = C3 x C3 x C3} \{1=7776, 2=28512\} \newline  [[0, 1, 2, 3], [0, 5, 18, 22], [0, 21, 23, $\infty$]]
\end{enumerate}
\end{example}

\section{1-rotational unitals of order 4}
For unitals of order $4$ together with "multiplier" filtering, calculation of automorphism group order and filtering by isomorphism was performed. According to this, $676$ non-isomorphic designs with automorphism group of order $64$ were found, which is bigger than $67$ in \cite{Krc}. For groups of order $128$, $192$, $256$, $384$, $768$ number of non-isomorphic designs is $75$, $5$, $12$, $1$, $3$, respectively. We didn't compare this data with \cite{NagyWeb} due to large amount of computation needed.
\begin{center}
\begin{tabular}{||c c c c||} 
 \hline
 GAP \# & Structure & \# of diff families & \# of non-iso designs \\ [0.5ex] 
 \hline\hline
 1 & {\tt C64} & 4 & 4 \\
 2 & {\tt C8 x C8} & 9 & 9 \\
 3 & {\tt C8 : C8} & 16 & 16 \\
 7 & {\tt Q8 : C8} & 93 & 68 \\
 11 & {\tt (C4 x C2) . (C4 x C2)} & 41 & 30 \\
 13 & {\tt (C4 x C2) . (C4 x C2)} & 42 & 42 \\
 14 & {\tt (C4 x C2) . (C4 x C2)} & 32 & 32 \\
 15 & {\tt C8 : C8} & 34 & 34 \\
 16 & {\tt C8 : C8} & 16 & 16 \\
 27 & {\tt C16 : C4} & 4 & 4 \\
 28 & {\tt C16 : C4} & 29 & 28 \\
 37 & {\tt (C4 x C2) . (C4 x C2)} & 27 & 26 \\
 43 & {\tt C8 . (C4 x C2)} & 2 & 2 \\
 44 & {\tt C4 : C16} & 4 & 4 \\
 45 & {\tt C8 . D8} & 6 & 6 \\
 46 & {\tt C16 : C4} & 15 & 15 \\
 48 & {\tt C16 : C4} & 6 & 6 \\
 49 & {\tt C4 . D16} & 3 & 3 \\
 120 & {\tt C4 x Q16} & 29 & 29 \\
 122 & {\tt Q16 : C4} & 27 & 27 \\
 126 & {\tt C8 x Q8} & 15 & 15 \\
 127 & {\tt C8 : Q8} & 40 & 40 \\
 143 & {\tt C4 : Q16} & 18 & 18 \\
 156 & {\tt Q8 : Q8} & 54 & 54 \\
 158 & {\tt Q8 : Q8} & 62 & 62 \\
 160 & {\tt (C2 x C2) . (C2 x D8)} & 67 & 67 \\
 168 & {\tt (C2 x C2) . (C2 x D8)} & 52 & 45 \\
 172 & {\tt (C2 x C2) . (C2 x D8)} & 24 & 24 \\
 175 & {\tt C4 : Q16} & 4 & 4 \\
 179 & {\tt C8 : Q8} & 24 & 23 \\
 180 & {\tt (C2 x C2) . (C2 x D8)} & 30 & 30 \\
 181 & {\tt C8 : Q8} & 12 & 12 \\
 182 & {\tt C8 : Q8} & 50 & 50 \\
 238 & {\tt Q8 : Q8} & 78 & 78 \\
 239 & {\tt Q8 x Q8} & 23 & 23 \\
 245 & {\tt (C2 x C2) . (C2 x C2 x C2 x C2)} & 21 & 21 \\
 \hline
\end{tabular}
\end{center}

\begin{example} The classical unital $S(2,5,65)$ appears as a 1-rotational design for groups {\tt SmallGroup(64,11)}, {\tt SmallGroup(64,28)} and {\tt SmallGroup(64,245)}.
\end{example}

\section{Enumerations of 1-rotational designs for k = 5}
1-rotational unitals of order $4$ $S(2,5,65)$ are listed in previous section. In this section we enumerate 1-rotational designs $S(2,k,v)$ with $k=5$ and $v\in \{25,85,101\}$, generated by actions of groups of order $24$, $84$ and $100$.
\begin{center}
\begin{tabular}{||c c c c||} 
 \hline
 GAP ID & Structure & \# of diff families & Comments \\ [0.5ex] 
 \hline\hline
 {\tt SmallGroup(25,1)} & {\tt C3 : C8} & 1 & - \\
 {\tt SmallGroup(25,2)} & {\tt C24} & 1 & - \\
 {\tt SmallGroup(25,3)} & {\tt SL(2,3)} & 1 & - \\
 \hline
 {\tt SmallGroup(84,1)} & {\tt C7 : C12} & 5626 & - \\
 {\tt SmallGroup(84,2)} & {\tt C4 x (C7 : C3)} & 689 & - \\
 {\tt SmallGroup(84,3)} & {\tt C7 x (C3 : C4)} & 1459 & - \\
 {\tt SmallGroup(84,4)} & {\tt C3 x (C7 : C4)} & 1673 & - \\
 {\tt SmallGroup(84,5)} & {\tt C21 : C4} & 1600 & - \\
 {\tt SmallGroup(84,6)} & {\tt C84} & 1072 & - \\
 {\tt SmallGroup(84,9)} & {\tt C2 x C2 x (C7 : C3)} & 503 & - \\
 {\tt SmallGroup(84,10)} & {\tt C7 x A4} & 180 & - \\
 {\tt SmallGroup(84,11)} & {\tt (C14 x C2) : C3} & 157 & - \\
 {\tt SmallGroup(84,15)} & {\tt C42 x C2} & 464 & - \\
 \hline
 {\tt SmallGroup(100,6)} & {\tt C5 x (C5 : C4)} & 110 & - \\
 {\tt SmallGroup(100,8)} & {\tt C20 x C5} & 40 & - \\
 {\tt SmallGroup(100,16)} & {\tt C10 x C10} & 98 & - \\
 \hline
\end{tabular}
\end{center}

\begin{example} There are no 1-rotational designs for $(v,k)\in\{(21,5),(41,5),(45,5),(61,5),(81,5)\}$.
\end{example}

\begin{example} Affine plane $(v,k)=(25,5)$ of order $5$ can be obtained in 3 ways.
\begin{enumerate}
\item {\tt SmallGroup(24,1) = C3 : C8} \{1=36000\} [[0, 1, 4, 14, 17], [0, 2, 3, 8, $\infty$]]
\item {\tt SmallGroup(24,2) = C24} \{1=36000\} [[0, 1, 2, 14, 22], [0, 3, 4, 11, $\infty$]]
\item {\tt SmallGroup(24,3) = SL(2,3)} \{1=36000\} [[0, 1, 6, 20, 21], [0, 3, 4, 11, $\infty$]]
\end{enumerate}
\end{example}

\section{Results location}
Raw results can be found at \href{https://github.com/Ihromant/math-utils/tree/master/src/test/resources/1rot}{https://github.com/Ihromant/math-utils/tree/master/src/test/resources/1rot}.
How to understand tables. If group is denoted with {\tt computer font} then groups were transformed to Cayley table using \cite{GAP4}, else they are in most obvious form (usually together with multipliers). Raw results are filtered by equivalence relation defined in Definition 1.1. When the number of an element equals the order of the group - this means that it is the fixed element of the group action.

\end{document}